\documentclass[12pt]{amsart}
\usepackage{calligra}
\usepackage{amsthm}
\usepackage{amsmath}
\usepackage{amssymb}
\usepackage{amscd}
\usepackage{bbm}
\usepackage[all]{xy}
\usepackage{mathrsfs}
\usepackage[top=72pt,bottom=96pt,left=72pt,right=72pt]{geometry}
\newtheorem{Lemma}{Lemma}[section]

\newtheorem{Definition}[Lemma]{Definition}

\def\Ad{\mathrm{Ad}}

\def\aut{\mathfrak{aut}}
\def\Aut{\mathrm{Aut}}

\def\End{\mathrm{End}}
\def\F{\mathscr{F}}

\def\g{\mathfrak{g}}

\def\hor{\mathrm{hor}}
\def\id{\mathrm{id}}
\def\im{\mathrm{im}}

\def\loc{\mathrm{loc}}

\def\P{\mathbb{P}}

\def\pfill{\par\vskip6pt plus4pt minus3pt\noindent}
\def\pr{\mathrm{pr}}

\def\R{\mathbb{R}}

\def\triv{\mathrm{triv}}
\def\vary#1{\delta\!#1}

\def\Vert{\mathrm{Vert}}
\def\vtriv{\mathrm{vtriv}}

\begin{document}
\date{\today}
\title{On Connections and their Curvatures}
\author{Gustavo Amilcar Salda\~na Moncada \& Gregor Weingart}
\address{Gustavo Amilcar Salda\~na Moncada\\
Instituto de Matem\'aticas (Ciudad de M\'exico) UNAM}
\email{gamilcar@ciencias.unam.mx}
\address{Gregor Weingart\\
Instituto de Matem\'aticas (Cuernavaca) UNAM}
\email{gw@matcuer.unam.mx}
\begin{abstract}
This paper presents a brief study on connections on fiber, principal and vector smooth bundles as well as some relations with their curvatures.
 \begin{center}
  \parbox{300pt}{\textit{MSC 2010:}\ 57R22.}
  \\[5pt]
  \parbox{300pt}{\textit{Keywords:}\ Fiber bundles, non--linear connections, curvature}
 \end{center}
\end{abstract}
\maketitle
\section{Introduction}

The study of bundles and their {\it ad hoc } connections is one of the must important subjects in differential geometry.
For example, the concept of vector bundle and linear connection  appear on every introductory course in the form of the tangent bundle and the Levi--Civita connection, respectively. Principal bundles also play an important role;  for example the frame bundle of a Riemannian manifold reflects different ways to give an orthonormal basis for tangent spaces, not to mentiontion that field theories in modern physics are developed by using principal bundles and principal connections \cite{gtvp}. All these objects arise from the more general concepts of fiber bundles and non--linear connections. This work's aim is to present a brief study on connections on fiber, principal and vector smooth bundles as well as some relations with their curvatures.\\

This paper corresponds to an extended version of the second section of \cite{ag}, but focused in geometrical relations, not the functioral ones and it breaks down as follows: in the second section we present the basics of the theory of fiber bundles and non--linear connections, as well as the relations of the curvature. The third section is about principal connections on principal bundles; while the fourth section is about vector bundles and linear connections. Finally in the appendix A we present some calculations related to the curvature and its relation with covariant derivatives in local coordinates.

\section{Non--Linear Connection}

A fiber bundle over a manifold $M$ with model fiber manifold
 $\F$ is a manifold $\F M$ endowed with a smooth projection map $\pi:\,
 \F M\longrightarrow M$, which is locally trivializable in the sense that
 there exists an open covering of $M$ by open subsets $U\,\subset\,M$,
 each of which allows for a diffeomorphism
 $$
  \Phi:\;\;U\;\times\;\F\;\stackrel\cong\longrightarrow\;\pi^{-1}(\,U\,)
 $$
 such that $\pi\,\circ\,\Phi:\,U\times\F\longrightarrow U$ equals the
 projection to $U$. The preimage of a point $p\,\in\,M$ under $\pi$ is called
 the fiber of the bundle over $p$, it is a submanifold $\F_pM\,:=\,\pi^{-1}
 (\,p\,)\,\subset\,\F M$ of the total space $\F M$ diffeomorphic to the model
 fiber $\F$. Morphisms between two fiber bundles are smooth maps $\varphi:\,\F M\longrightarrow
 \hat\F M$ between the total spaces which commute with the respective
 projections
 \begin{equation}\label{fbmor}
  \vcenter{\hbox{\begin{picture}(90,50)(0,0)
   \put( 0,40){$\F M$}
   \put(68,40){$\hat\F M$}
   \put(40, 0){$M$}
   \put(26,43){\vector(+1, 0){39}}
   \put(43,47){$\scriptstyle\varphi$}
   \put(15,37){\vector(+1,-1){26}}
   \put(20,21){$\scriptstyle\pi$}
   \put(78,37){\vector(-1,-1){26}}
   \put(69,21){$\scriptstyle\hat\pi$}
  \end{picture}}}
 \end{equation}
 and thus map the fibers of $\F M$ to the fibers of $\hat\F M$ over the same
 point.\\  
 
 Intuitively, a non--linear or Ehresmann
 connection on a fiber bundle $\F M$ is a left hand side section $\nabla$
 of the short exact sequence of vector bundles
 \begin{equation}\label{nlc}
  0\;\longrightarrow\;
  \hbox to0pt{\begin{picture}(112,27)(0,0)
   \qbezier(89,11)(66,22)(45,11)
   \put(45,11){\vector(-2,-1){0}}
   \put(60,20){$\nabla$}
  \end{picture}\hss}
  \Vert\;\F M
  \;\stackrel\subset\longrightarrow\;
  T\,\F M
  \;\stackrel{\pi_*}\longrightarrow\;
  \pi^*TM
  \;\longrightarrow\;0
 \end{equation}
 over $\F M$, the corresponding right hand side section is appropriately
 called the horizontal lift associated to $\nabla$. In order to studying connections
 we prefer the following version of this intuitive notion of
 a non--linear connection:

 \begin{Definition}[Non--linear Connections on Fiber Bundles]
 \hfill\label{cfb}\break
  A general or a non--linear connection on a fiber bundle $\F M$ over a manifold $M$
  is a field $\P^\nabla$ $\in$ $\Gamma(\,\F M,\,\End\,T\F M\,)$ of projections
  $(\P^\nabla)^2\,=\,\P^\nabla$ on the tangent bundle $T\F M$ such that its
  image distribution equals the vertical foliation:
  $$
   \im\;\Big(\;\P^\nabla_f:\;\;T_f\F M\;\longrightarrow\;T_f\F M\;\Big)
   \;\;\stackrel!=\;\;
   \Vert_f\F M\ .
  $$
 \end{Definition}

 Every non--linear connection $\P^\nabla$ on a fiber bundle $\F M$ allows us
 to define the directional derivative $D^\nabla_Xf\,\in\,\Gamma_\loc(\,M,\,
 \Vert\,\F M\,)$ of a given local section $f\,\in\,\Gamma_\loc(\,M,\,\F M\,)$
 in the direction of a vector field $X$ on $M$ by:
 $$
  (\;D^\nabla_Xf\;)_p
  \;\;:=\;\;
  \Big(\;T_pM\;\stackrel{f_{*,\,p}}\longrightarrow\;T_{f(p)}\F M
  \;\stackrel{\P^\nabla_{f(p)}}\longrightarrow\;\Vert_{f(p)}\F M\;\Big)
  \;X_p\ .
 $$
 These directional derivatives assemble into a first order differential
 operator
 \begin{equation}\label{nabla}
  D^\nabla:\;\;\Gamma(\,M,\,TM\,)\;\times\;\Gamma_\loc(\,M,\,\F M\,)
  \;\longrightarrow\;\Gamma_\loc(\,M,\,\Vert\,\F M\,)\ ,
 \end{equation}
 which is the non--linear analogue of the classical definition of covariant
 derivatives on vector bundles. Somewhat annoyingly this covariant derivative
 $D^\nabla_Xf$ contains the redundant information $f\,=\,\pi_{\F M}\,\circ\,
 D^\nabla_Xf$, where $\pi_{\F M}$ denotes the vertical tangent bundle
 projection $\Vert\,\F M\longrightarrow\F M$. The simplicity of linear
 and principal connections stems from the fact that we can get rid of
 this redundancy altogether, the reduced covariant derivative $\nabla_Xf$
 captures only the partial derivatives of the section $f$.\\

 The Nijenhuis or curvature tensor of a non--linear connection $\P^\nabla$
 on a fiber bundle $\F M$ over a manifold $M$ is the horizontal $2$--form
 $R^\nabla$ on the total space $\F M$ of the fiber bundle with values in
 the vertical tangent bundle defined for two arbitrary vector fields
 $X,\,Y$ on $\F M$ by:
 \begin{eqnarray*}
  \lefteqn{R^\nabla(\;X,\,Y\;)}
  &&
  \\[4pt]
  &:=&
  -\;\P^\nabla\,[\;(\,\id-\P^\nabla\,)\,X,\;(\,\id-\P^\nabla\,)\,Y\;]
  \;-\;(\,\id-\P^\nabla\,)\,[\;\P^\nabla X,\;\P^\nabla Y\;]
  \\[4pt]
  &=&
  -\,\P^\nabla\,[\,X,\,Y\,]
  \,+\,\P^\nabla\,[\,X,\,\P^\nabla Y\,]
  \,+\,\P^\nabla\,[\,\P^\nabla X,\,Y\,]
  \,-\,[\,\P^\nabla X,\P^\nabla Y\,]\ .
 \end{eqnarray*}
 The strange sign is necessary to make this definition agree with the
 classical definition of the curvature of a linear connection, compare
 equation (\ref{cross}). By construction $R^\nabla\,\in\,\Omega^2_\hor
 (\,\F M,\,\Vert\,\F M\,)$ is $C^\infty(\,\F M\,)$--bilinear and alternating.
 The integrability constraint for the vertical foliation $\Vert\,\F M$ tells
 us $(\,\id\,-\,\P^\nabla\,)\,[\,\P^\nabla X,\,\P^\nabla Y\,]\,=\,0$,
 and the simplified expression
 \begin{equation}\label{rcurv}
  R^\nabla(\;X,\;Y\;)
  \;\;=\;\;
  -\;\P^\nabla\,
  [\;(\,\id\,-\,\P^\nabla\,)\,X,\;(\,\id\,-\,\P^\nabla\,)\,Y\;]\ ,
 \end{equation}
 evidently results in a vertical vector field for all vector fields $X,\,Y$
 on $\F M$ and vanishes for a vertical vector field $X$ or $Y$ no matter what
 is the other. In particular the curvature $R^\nabla$ measures exactly the
 failure of the horizontal distribution $\ker\,\P^\nabla\,\subseteq\,T\F M$
 associated to $\P^\nabla$ to be integrable.\\
 
 An interpretation of the curvature tensor along classical lines as a
 commutator of covariant derivatives requires a significant amount of
 extra work and so we will only sketch the basic idea. The iterated tangent
 bundle of a smooth manifold $\F$ is a smooth manifold $T(\,T\F\,)$ coming
 along with a canonical involutive diffeomorphism $\Theta:\,T(\,T\F\,)
 \longrightarrow T(\,T\F\,)$ characterized by
 \begin{equation}\label{theta}
  \Theta\;\left(\;\left.\frac d{dt}\right|_0\;
  \left.\frac d{d\varepsilon}\right|_0\gamma(\,t,\,\varepsilon\,)\;\right)
  \;\;=\;\;
  \left(\;\left.\frac d{dt}\right|_0\;\left.\frac d{d\varepsilon}\right|_0
  \gamma(\,\varepsilon,\,t\,)\;\right)
 \end{equation}
 for all $\gamma:\,\R^2\longrightarrow\F$. By construction, $\Theta$ fits
 into the commutative diagram
 \begin{equation}\label{swap}
  \vcenter{\hbox{\begin{picture}(170,62)(0,0)
   \put( 10,52){$T(\,T\F\,)$}
   \put(125,52){$T(\,T\F\,)$}
   \put(  0, 0){$T\F\oplus T\F$}
   \put(115, 0){$T\F\oplus T\F$}
   \put( 57,56){\vector(+1, 0){52}}
   \put( 78,60){$\scriptstyle\Theta$}
   \put( 57, 4){\vector(+1, 0){52}}
   \put( 74, 9){$\scriptstyle\mathrm{swap}$}
   \put( 28,48){\vector( 0,-1){38}}
   \put( 18,28){$\scriptstyle\Pi$}
   \put(143,48){\vector( 0,-1){38}}
   \put(147,28){$\scriptstyle\Pi$}
  \end{picture}}}
 \end{equation}
 where $\mathrm{swap}$ interchanges the two summands and $\Pi$ is the
 double projection:
 $$
  \Pi\left(\;\left.\frac d{dt}\right|_0\;\left.\frac d{d\varepsilon}\right|_0
  \gamma(\,t,\varepsilon\,)\;\right)
  \;\;:=\;\;
  \left.\frac d{d\varepsilon}\right|_0\gamma(\,0,\,\varepsilon\,)
  \;\oplus\;\left.\frac d{dt}\right|_0\gamma(\,t,\,0\,)\ .
 $$
 The involution $\Theta$ and double projection $\Pi$ generalize
 directly to iterated vertical tangent bundles, because we may identify
 the fibers of $\Vert\,\Vert\,\F M$ considered as a fiber bundle over
 $M$ in every point $p\,\in\,M$ with the iterated tangent bundle
 $T(\,T[\,\F_pM\,]\,)\,=\,[\,\Vert\,\Vert\,\F M\,]_p$. The involutions and
 double projections thus defined on all the fibers of $\Vert\,\Vert\,\F M$
 assemble into an involutive fiber bundle automorphism $\Theta$ and a double
 projection $\Pi:\,\Vert\,\Vert\,\F M\longrightarrow\Vert\,\F M\,\oplus\,
 \Vert\,\F M$. Making good use of the involution $\Theta$ we extend a
 non--linear connection $\P^\nabla$ on $\F M$ to a unique non--linear
 connection $\P^{\nabla^\Vert}$ on $\Vert\,\F M$ by stipulating the identity
 \begin{equation}\label{vconex}
  D^{\nabla^\Vert}_X
  \Big(\;\left.\frac d{d\varepsilon}\right|_0f_\varepsilon\;\Big)
  \;\;=\;\;
  \Theta\,
  \Big(\;\left.\frac d{d\varepsilon}\right|_0D^\nabla_Xf_\varepsilon\;\Big)
 \end{equation}
 for every vector field $X\,\in\,\Gamma(\,M,\,TM\,)$ and every
 smooth one--parameter family $(\,f_\varepsilon\,)_{\varepsilon\,\in\,\R}$
 of local sections $f_\varepsilon$ of $\F M$ with infinitesimal variation
 $\left.\frac d{d\varepsilon}\right|_0f_\varepsilon\,\in\,\Gamma_\loc
 (\,M,\,\Vert\,\F M\,)$. In this construction of the induced connection
 on $\Vert\,\F M$, the involution $\Theta$ is proper remedy for the nuisance
 $$
  \begin{array}{ccccccl}
   \Pi\,\Big(\!\!\!&D^{\nabla^\Vert}_X
   (\,\left.\frac d{d\varepsilon}\right|_0f_\varepsilon\,)&\!\!\!\Big)
   \;\;=\;\;
   \Big(\!\!\!&\left.\frac d{d\varepsilon}\right|_0f_\varepsilon&\!\!\!\Big)
   \,\oplus\,\Big(\!\!\!&D^\nabla_Xf_0&\!\!\!\Big)
   \\[7pt]
   \Pi\,\Big(\!\!\!&\left.\frac d{d\varepsilon}\right|_0
   D^\nabla_Xf_\varepsilon&\!\!\!\Big)
   \;\;=\;\;
   \Big(\!\!\!&D^\nabla_Xf_0&\!\!\!\Big)\,\oplus\,\Big(\!\!\!&
   \left.\frac d{d\varepsilon}\right|_0f_\varepsilon&\!\!\!\Big)\ ,
  \end{array}
 $$
 compare the commutative diagram (\ref{swap}). For essentially the same
 reason the involution $\Theta$ appears in the key identity linking the
 curvature $R^\nabla$ of a non--linear connection $\P^\nabla$ on a fiber
 bundle $\F M$ over $M$ to the commutator
 \begin{equation}\label{cross}
  R^\nabla_{X,\,Y}f
  \;\;=\;\;
  \Big(\;D^{\nabla^\Vert}_XD^{\nabla^{\hphantom{V}}}_Y\!\!f
  \;-\;\Theta(\;D^{\nabla^\Vert}_YD^{\nabla^{\hphantom{V}}}_X\!\!f\;)\;\Big)
  \;-\;D^{\nabla^{\hphantom{V}}}_{[\,X,\,Y\,]}f
 \end{equation}
 of covariant derivatives of a local section $f\,\in\,\Gamma_\loc(\,M,\,
 \F M\,)$ in the direction of vector fields $X,\,Y\,\in\,\Gamma(\,M,\,TM\,)$,
 in which $R^\nabla_{X,\,Y}f$ is a simplified notation for the local section
 of $\Vert\,\F M$ over $M$ defined in $p\,\in\,M$ by:
 $$
  \Big(\;R^\nabla_{X,\,Y}f\;\Big)_p
  \;\;:=\;\;
  R^\nabla_{f(p)}(\;f_{*,\,p}X_p,\;f_{*,\,p}Y_p\;)
  \;\;\in\;\;
  \Vert_{f(p)}\F M\ .
 $$
 The basic tenet of symbolic calculus on jet bundles provides the correct
 interpretation for the differences in (\ref{cross}): The fibers of the
 double projection $\Pi:\,\Vert\,\Vert\,\F M\longrightarrow\Vert\,\F M\oplus
 \Vert\,\F M$ are affine spaces modelled on the fibers of $\Vert\,\F M$ and
 thus come along with a difference map
 $$
  -\;:\;\;
  \Vert\,\Vert\,\F M\;\times_{\Vert\,\F M\oplus\Vert\,\F M}\Vert\,\Vert\,\F M
  \;\longrightarrow\;\Vert\,\F M\ ,
 $$
 which determines the inner difference in the curvature identity (\ref{cross}),
 the outer difference is simply the difference in $\Vert\,\F M$ considered as
 a vector bundle over $\F M$. Taking care with this delicate interpretation of
 differences we will prove the curvature identity (\ref{cross}) in local
 coordinates in Appendix \ref{excalc}.
 
 \begin{Definition}[Parallel morphisms between Fiber Bundles]
 \hfill\label{ph}\break
  A parallel morphism between fiber bundles $\F M$ and $\hat\F M$ over the
  same manifold $M$ endowed with connections $\P^\nabla$ and $\P^{\hat\nabla}$
  respectively is a morphism $\varphi:\,\F M\longrightarrow\hat\F M$ of
  fiber bundles  such that the following
  diagram commutes:
  $$
   \begin{picture}(116,55)(0,0)
    \put(  0,46){$T\,\F M$}
    \put( 81,46){$T\,\hat\F M$}
    \put(  0, 0){$T\,\F M$}
    \put( 81, 0){$T\,\hat\F M$}
    \put( 37,50){\vector(+1, 0){42}}
    \put( 55,54){$\scriptstyle\varphi_*$}
    \put( 37, 4){\vector(+1, 0){42}}
    \put( 55, 8){$\scriptstyle\varphi_*$}
    \put( 16,42){\vector( 0,-1){30}}
    \put(  1,25){$\scriptstyle\P^\nabla$}
    \put( 96,42){\vector( 0,-1){30}}
    \put(100,25){$\scriptstyle\P^{\hat\nabla}$}
   \end{picture}
  $$
 \end{Definition}

 \noindent
 The constraint $\hat\pi\,\circ\,\varphi\,=\,\pi$ characterizing morphisms
 of fiber bundles  readily implies $$\varphi_*(\,\Vert
 \,\F M\,)\,\subset\,\Vert\,\hat\F M,$$ hence a morphism $\varphi$ of
 fiber bundles is parallel, if and only if $\varphi_*$ maps the horizontal
 distribution of $\F M$ to the horizontal distribution of $\hat\F M$:
 $$
  \varphi\;\textrm{\ parallel}
  \qquad\Longleftrightarrow\qquad
  \varphi_*(\;\ker\;\P^\nabla\;)
  \;\;\subset\;\;
  \ker\;\P^{\hat\nabla}\ .
 $$

\section{Principal Connections}

 Having discussed general non--linear connections on fiber bundles in some
 detail we now  want to specialize to principal  connections on principal bundles in this section.  Recall first of all that a principal bundle
 modelled on a Lie group $G$ is a smooth fiber bundle $GM$ with model fiber
 $G$ endowed with a smooth right, fiber preserving action of $G$ on its
 total space $GM$
 $$
  \rho\,:\;\;GM\;\times\;G\;\longrightarrow\;GM,
  \qquad(\,g\,,\,\gamma\,)\;\longmapsto\;g\,\gamma\ ,
 $$
 which is simply transitive on all fibers of $GM$ over $M$. Every local
 section ${\scriptstyle\Gamma}:\,U\longrightarrow GM$ of the projection
 $\pi:\,GM\longrightarrow M$ defined over an open subset $U$ extends to
 a $G$--equivariant local trivialization of $GM$ over $U$
 \begin{equation}\label{pbtriv}
  \Phi:\;\;U\times G
  \;\stackrel{{\scriptscriptstyle\Gamma}\times\id}\longrightarrow\;
  \pi^{-1}(\,U\,)\times G
  \;\stackrel\rho\longrightarrow\;
  \pi^{-1}(\,U\,)\ ,
 \end{equation}
 satisfying $\Phi(u,g\gamma)\,=\,{\scriptstyle\Gamma}(u)g\gamma\,=\,
 \Phi(u,g)\,\gamma$ by construction. The existence of $G$--equivariant
 trivializations ensures that the unique set theoretic map
 \begin{equation}\label{bs}
  \backslash\;:\;\;GM\;\times_M\,GM\;\longrightarrow\;G,
  \qquad(\,g,\,\hat g\,)\;\longmapsto\;g^{-1}\,\hat g\ ,
 \end{equation}
 satisfying $g\,(g^{-1}\hat g)\,=\,\hat g$ for all $g,\,\hat g\,\in\,GM$ in
 the same fiber is actually a smooth map, because it reads $$\Phi(u,g)^{-1}
 \Phi(u,\hat g)\,=\,g^{-1}\hat g$$ in every such trivialization $\Phi$.\\
 
  The
 automorphism group bundle of a principal bundle $GM$ over a manifold $M$
 is the Lie group bundle $\Aut\;GM$ over $M$ defined by
 \begin{equation}\label{aut}
  \Aut\,GM
  \;\;:=\;\;
  \{\;\;(p,\psi)\;\;|\;\;\psi:\,G_pM\longrightarrow G_pM
  \textrm{\ is $G$--equivariant}\;\;\}
 \end{equation}
 with the bundle projection $\pi_{\Aut\,GM}\,:\,\Aut\,GM\longrightarrow M,
 \,(\,p,\,\psi\,)\longmapsto p$. In mathematical physics, the Fr\'echet--Lie
 group $\Gamma(\,M,\,\Aut\,GM\,)$ of all global sections of the automorphism
 bundle is called the gauge group of $GM$.\\

 The fiber of the Lie group bundle $\Aut\;GM$ over a point $p\,\in\,M$ is a
 Lie group $\Aut_pGM$ isomorphic, although not canonically so, to the original
 group $G$, in particular its Lie algebra $\aut_pGM\,\cong\,\g$ is
 isomorphic to the Lie algebra of $G$. All these Lie algebras assemble
 into a smooth Lie algebra bundle $\aut\,GM$, whose global sections
 $\Gamma(\,M,\,\aut\,GM\,)$ form the Fr\'echet--Lie algebra of the gauge
 group $\Gamma(\,M,\,\Aut\,GM\,)$ of the principal bundle $GM$.

 \begin{Definition}[Principal Connections]
 \hfill\label{pconex}\break
  A principal connection on a principal $G$--bundle $GM$ over a manifold
  $M$ is a non--linear connection $\P^\nabla$ on the fiber bundle $GM$,
  which is invariant under the right action of $G$ on $GM$ in the sense
  that the right translations $R_\gamma:\,GM\longrightarrow GM,\,g
  \longmapsto g\gamma,$ are parallel automorphisms for all $\gamma\,\in\,G$.
 \end{Definition}
 
 \noindent
 In difference to general fiber bundles, the vertical tangent bundle
 of a principal bundle $GM$ is trivializable: The construction of the
 Maurer--Cartan form on a Lie group $G$ applies verbatim and results
 in the vertical trivialization isomorphism $$\Vert\,GM\stackrel\cong
 \longrightarrow GM\times\g,\;\left.\frac d{dt}\right|_0g_t\longmapsto
 (g_0,\left.\frac d{dt}\right|_0g^{-1}_0g_t)$$ under the proviso that
 the curve $t\longmapsto g_t$ chosen to represent the vertical tangent
 vector stays in the fiber so that the curve $t\longmapsto g_0^{-1}g_t$
 is defined by application (\ref{bs}). The composition of this vertical
 trivialization with the projection to $\g$
 \begin{equation}\label{pco}
  \vcenter{\hbox{\begin{picture}(230,54)(0,0)
   \put(  0,40){$\vtriv\,:$}
   \put( 40,40){$\Vert\,GM$}
   \put(133,40){$GM\times\g$}
   \put( 95, 0){$GM$}
   \put(226,40){$\g$}
   \put( 89,43){\vector(+1, 0){40}}
   \put(105,46){$\scriptstyle\cong$}
   \put( 66,36){\vector(+1,-1){27}}
   \put( 56,21){$\scriptstyle\pi_{GM}$}
   \put(146,36){\vector(-1,-1){27}}
   \put(139,21){$\scriptstyle\pr_L$}
   \put(182,43){\vector(+1, 0){40}}
   \put(194,48){$\scriptstyle\pr_R$}
  \end{picture}}}
 \end{equation}
 induces isomorphisms $\vtriv_g:\,\Vert_gGM\stackrel\cong\longrightarrow\g$
 at every $g\,\in\,GM$ and thus a bijection $\P^\nabla\,\longleftrightarrow\,
 \omega$ between non--linear connections and $\g$--valued $1$--forms
 $\omega\,\in\,\Omega^1(\,GM,\,\g\,)$, which agree with $\vtriv$ on
 vertical tangent vectors:
 \begin{equation}\label{omega}
  \omega
  \;\;:=\;\;
  \vtriv\,\circ\,\P^\nabla
  \qquad\Longleftrightarrow\qquad
  \P^\nabla
  \;\;:=\;\;
  \vtriv^{-1}\,\circ\,\omega\ .
 \end{equation}
 This description of non--linear connections on $GM$ by $\g$--valued
 $1$--forms $\omega$ with $\left.\omega\right|_{\Vert\,GM}\,=\,\vtriv$
 is much more convenient and will be used for the rest of this article.\\
 
 To characterize the principal connections specified in Definition
 \ref{pconex} by their connection forms $\omega$ we recall that the
 right translation by $\gamma\,\in\,G$ is a parallel morphism
 $$R_\gamma:\,GM\longrightarrow GM,\;g\longmapsto g\gamma,$$ if and only
 if its differential $R_{\gamma*}$ maps the horizontal vectors
 $\ker\,\omega_g\,\subset\,T_gGM$ at a point $g\,\in\,GM$ to horizontal
 vectors $\ker\omega_{g\gamma}\,\subset\,T_{g\gamma}GM$. On the
 complementary vertical vectors the connection form behaves like
 $\vtriv\,=\,\left.\omega\right|_{\Vert\,GM}$:
 $$
  (\vtriv\,\circ\,R_{\gamma*})\Big(\!\left.\frac d{dt}\right|_0g_t\Big)
  \;\;=\;\;
  \left.\frac d{dt}\right|_0\!(\gamma g_0)^{-1}(g_t\gamma)
  \;\;=\;\;
  \Ad_{\gamma^{-1}}\vtriv\Big(\!\left.\frac d{dt}\right|_0g_t\Big)\ .
 $$
 In consequence the condition $\omega\,\circ\,R_{\gamma*}\,=\,
 \Ad_{\gamma^{-1}}\,\circ\,\omega$ for all $\gamma\,\in\,G$ is
 necessary and sufficient for a $\g$--valued $1$--form $\omega$ agreeing
 on vertical vectors with $\vtriv\,=\,\left.\omega\right|_{\Vert\,GM}$
 to be the connection form of a principal connection $\P^\nabla$. The
 chain rule in the guise of a Leibniz rule tells us for such a form
 \begin{eqnarray*}
  \omega_{g_0\gamma_0}\Big(\;\left.\frac d{dt}\right|_0g_t\,\gamma_t\;\Big)
  &=&
  \omega_{g_0\gamma_0}\Big(\;\left.\frac d{dt}\right|_0g_t\,\gamma_0
  \;+\;\left.\frac d{dt}\right|_0g_0\,\gamma_t\;\Big)
  \\
  &=&
  \omega_{g_0\gamma_0}\Big(\;(\,R_{\gamma_0}\,)_*\left.\frac d{dt}
  \right|_0g_t\;\Big)\;+\;\left.\frac d{dt}\right|_0(g_0\gamma_0)^{-1}
  (g_0\gamma_t)
  \\
  &=&
  \Ad_{\gamma^{-1}_0}\;\omega_{g_0}\Big(\;\left.\frac d{dt}\right|_0g_t\;\Big)
  \;+\;\left.\frac d{dt}\right|_0\gamma_0^{-1}\gamma_t
 \end{eqnarray*}
 for every curve $t\longmapsto g_t$ in $GM$ and $t\longmapsto\gamma_t$ in $G$.
 Put differently:

 \begin{Lemma}[Principal Connection Axiom]
 \hfill\label{cpb}\break
  On every principal $G$--bundle $GM$ the association $\P^\nabla\,
  \longleftrightarrow\,\omega$ characterized by $\omega\,:=\,\vtriv\circ
  \P^\nabla$ induces a bijection between principal connections in the sense
  of Definition \ref{pconex} and $\g$--valued $1$--forms $\omega$ on $GM$
  satisfying the axiom
  $$
   \omega_{g_0\gamma_0}\Big(\;\left.\frac d{dt}\right|_0g_t\,\gamma_t\;\Big)
   \;\;=\;\;
   \Ad_{\gamma^{-1}_0}\;\omega_{g_0}\Big(\;\left.\frac d{dt}\right|_0g_t\;\Big)
   \;+\;\left.\frac d{dt}\right|_0\gamma_0^{-1}\gamma_t
  $$
  for all choices of smooth curves $t\longmapsto g_t$ in $GM$ and curves
  $t\longmapsto\gamma_t$ in $G$.
 \end{Lemma}

 \noindent
 Cartan's Second Structure Equation \cite{gtvp} is a convenient description
 of the image of the composition of the curvature tensor $R^\nabla$ with the
 vertical trivialization $\vtriv$ in terms of the exterior derivative of the
 connection form
 \begin{equation}\label{cstruc}
  \Omega
  \;\;:=\;\;
  \vtriv\,\circ\,R^\nabla
  \;\;\stackrel!=\;\;
  d\omega\;+\;\frac12\,[\,\omega\,\wedge\,\omega\,]\ ,
 \end{equation}
 where $\frac12[\,\omega\wedge\omega\,](X,Y)\,:=\,[\,\omega(X),\,\omega(Y)\,]$;
 the usefulness of the factor $\frac12$ in this definition is beyond question.
 Cartan's Second Structure Equation can be easily verified directly. On the
 other hand it can be derived without too much additional effort from the
 curvature identity (\ref{cross}) by considering the reduced covariant
 derivative associated to a principal connection $\omega$
 $$
  \nabla:\;\;\Gamma(\,M,\,TM\,)\,\times\,\Gamma_\loc(\,M,\,GM\,)
  \;\longrightarrow\;C^\infty_\loc(\,M,\,\g\,)
 $$
 defined by $\nabla_Xf\,:=\,\vtriv(\,D^\nabla_Xf\,)$; in light of equation
 (\ref{omega}) this becomes
 $$
  (\;\nabla_Xf\;)_p
  \;\;:=\;\;
  \vtriv(\;\P^\nabla_{f(p)}(\,f_{*,\,p}X_p\,)\;)
  \;\;\stackrel!=\;\;
  (\,f^*\omega\,)_p(\,X\,)
 $$
 for all $X\,\in\,\Gamma(M,TM)$. Iterating the vertical trivialization
 isomorphism $\Vert\,GM\,\longrightarrow GM\times\g$ used to define
 $\vtriv$ we obtain the isomorphism:
 \begin{equation}\label{itvtriv}
  \Psi:\;\Vert\,\Vert\,GM
  \;\stackrel\cong\longrightarrow\;
  \Vert(\,GM\times\g\,)
  \;\stackrel\cong\longrightarrow\;
  (GM\times\g)\,\times\,(\g\times\g)\ .
 \end{equation}
 Its inverse can be written down explicitly by unravelling the
 definitions
 \begin{eqnarray*}
  \Psi^{-1}
  (\,g,\,X,\,Y,\,Z\,)
  &=&
  \left.\frac d{dt}\right|_0\,\left.\frac d{d\varepsilon}\right|_0
  g\,e^{tX}\,e^{\varepsilon(\,Y\,+\,tZ\,)}
  \\
  &=&
  \left.\frac d{dt}\right|_0\,\left.\frac d{d\varepsilon}\right|_0
  g\,e^{tX\,+\,\varepsilon Y\,+\,t\varepsilon(\,Z\,+\,\frac12\,[X,Y]\,)}\ ,
 \end{eqnarray*}
 the second equality is the Baker--Campbell--Hausdorff Formula \cite{nodg}
 in the Lie group $G$ ignoring as usual terms of order $O(\,t^2,\,
 \varepsilon^2\,)$. In particular the involution $\Theta$ defined in equation
 (\ref{theta}) simply by the interchange $t\,\leftrightarrow\,\varepsilon$
 is not as innocent as it may appear, under the isomorphism $\Psi$ it becomes:
 \begin{equation}\label{picks}
  \Theta(\;g,\;X;\;Y,\;Z\;)
  \;\;\stackrel!=\;\;
  (\;g,\;Y;\;X,\;Z\,+\,[\,X,\,Y\,]\;)\ .
 \end{equation}
 In accordance with \cite{ag} the vertical trivialization isomorphism $\Vert\,GM\longrightarrow
 GM\times\g$ is parallel with respect to the connections $\P^{\nabla^\Vert}$
 and $\P^\nabla\oplus\P^\triv$  on $\Vert\,GM$ and $GM\times\g$ respectively,
 using the fiber bundle isomorphism $\Psi$ we thus obtain
 \begin{eqnarray*}
  \Psi(\;D^{\nabla^\Vert}_X\,D^{\nabla^{\hphantom{V}}}_Y\!\!f\;)
  &=&
  (\;f,\;(f^*\omega)(X);\;(f^*\omega)(Y),\;X\,(f^*\omega)(Y)\;)
  \\[2pt]
  \Psi(\;D^{\nabla^\Vert}_Y\,D^{\nabla^{\hphantom{V}}}_X\!\!f\;)
  &=&
  (\;f,\;(f^*\omega)(Y);\;(f^*\omega)(X),\;Y\,(f^*\omega)(X)\;)
 \end{eqnarray*}
 for all local sections $f\,\in\,\Gamma_\loc(\,M,\,GM\,)$ and all $X,\,Y
 \,\in\,\Gamma(\,M,\,TM\,)$. Evidently the double projections under $\Pi$
 of these two local sections of $\Vert\,\Vert\,GM$ do not agree unless we
 apply the involution $\Theta$ first picking up the rather unexpected
 additional term of equation (\ref{picks}) along the way:
 \begin{eqnarray*}
  \lefteqn{D^{\nabla^\Vert}_X\,D^{\nabla^{\hphantom{V}}}_Y\!\!f
   \;-\;\Theta(\;D^{\nabla^\Vert}_Y\,D^{\nabla^{\hphantom{V}}}_X\!\!f\;)}
  \quad
  &&
  \\
  &=&
  \Big(\;f,\;X(f^*\omega)(Y)\;-\;Y(f^*\omega)(X)\;-\;[\,(f^*\omega)(Y),\,
  (f^*\omega)(X)\,]\;\Big)
 \end{eqnarray*}
 Subtracting $D^\nabla_{[X,Y]}f\,=\,(\,f,\,(f^*\omega)([X,Y])\,)$ we obtain
 eventually
 \begin{eqnarray*}
  R^\nabla_{X,\,Y}f
  &=&
  (\;f,\;d(f^*\omega)(X,Y)\;+\;\frac12[\,(f^*\omega)\,\wedge\,(f^*\omega)\,]
  (\,X,\,Y\,)\;)
  \\
  &=&
  (\;f,\;(\,f^*\Omega\,)(\,X,\,Y\,)\;)
 \end{eqnarray*}
 for all local sections $f\,\in\,\Gamma_\loc(\,M,GM\,)$ and all $X,\,Y\,\in\,
 \Gamma(\,M,TM\,)$ by using the naturality $d(f^*\omega)\,=\,f^*(d\omega)$ of
 the exterior derivative $d$ and the definition of the $2$--form $\Omega$ in
 Cartan's Second Structure Equation (\ref{cstruc}).

\section{Linear Connections}

\pfill
 The strategy persued for linear connections on vector bundles $VM$ follows
 the model of principal connections closely. The tangent bundle of a vector
 space is canonically trivializable $TV\,\cong\,V\times V$ by taking
 actual derivatives and this becomes via $[\,\Vert\,VM\,]_p\,=\,T(\,V_pM\,)$
 the vertical trivialization $$\Vert\,VM\stackrel\cong\longrightarrow VM
 \oplus VM,\;\left.\frac d{dt}\right|_0v_t\longmapsto v_0\oplus\lim_{t\to0}
 \frac1t(v_t-v_0),$$ again under the proviso that the representative curve
 $t\longmapsto v_t$ is chosen to stay in the same vector space so that the
 quotient $\frac1t(v_t-v_0)\,\in\,V_{\pi(v_0)}M$ is well--defined. Composing
 with the projection to the right factor we obtain
 \begin{equation}\label{linco}
  \vcenter{\hbox{\begin{picture}(252,54)(0,0)
   \put(  0,40){$\vtriv:$}
   \put( 40,40){$\Vert\,VM$}
   \put(130,40){$VM\oplus VM$}
   \put(234,40){$VM$}
   \put( 97, 0){$VM$\ ,}
   \put( 88,43){\vector(+1, 0){40}}
   \put(105,46){$\scriptstyle\cong$}
   \put( 71,36){\vector(+1,-1){26}}
   \put( 60,21){$\scriptstyle\pi_{VM}$}
   \put(146,36){\vector(-1,-1){26}}
   \put(138,21){$\scriptstyle\pr_L$}
   \put(191,43){\vector(+1, 0){40}}
   \put(205,48){$\scriptstyle\pr_R$}
  \end{picture}}}
 \end{equation}
 which can be used to project out $\nabla_Xv\,:=\,\vtriv(D^\nabla_Xv)$ the
 redundant information from the covariant derivative $D^\nabla_Xv$ of a
 section $v\,\in\,\Gamma(\,M,\,VM\,)$:
  
 \begin{Definition}[Linear Connections on Vector Bundles]
 \hfill\label{linx}\break
  A linear connection on a vector bundle $VM$ on $M$ is a non--linear
  connection $\P^\nabla$ on $VM$ such that the reduced covariant derivative
  is $\R$--bilinear:
  $$
   \nabla:\;\;\Gamma(\,M,\,TM\,)\;\times\;\Gamma(\,M,\,VM\,)
   \;\longrightarrow\;\Gamma(\,M,\,VM\,)\ .
  $$
 \end{Definition}

 In order to characterize linear connections by a condition easier to verify
 we consider a vector bundle $VM$ with projection $\pi:\,VM\longrightarrow M$
 and a local basis of sections $\mathbf{v}_1,\,\ldots,\,\mathbf{v}_n\,\in\,
 \Gamma_\loc(\,M,\,VM\,)$ defined over the open domain $U\,\subseteq\,M$ of
 a local coordinate chart $x:\,U\longrightarrow x(U)\,\subseteq\,\R^m$. In
 this situation we can construct linear local coordinates for the total space
 $VM$ of the vector bundle $(x,v):\,\pi^{-1}(U)\longrightarrow x(U)\times\R^n$
 characterized by:
 $$
  (\,x^1,\,\ldots\,x^m;\,v^1,\,\ldots,\,v^n\,)
  \;\;\longleftrightarrow\;\;
  v^1\,\mathbf{v}_1\;+\;\ldots\;+\;v^n\,\mathbf{v}_n
  \;\;\in\;\;
  V_{(x^1,\ldots,x^m)}M\ .
 $$
 In these linear local coordinates a general non--linear connection $\P^\nabla$
 reads
 \begin{equation}\label{gnlc}
  \P^\nabla_{(x,v)}\Big(\frac\partial{\partial v^\alpha}\Big)
  \;\;=\;\;
  \frac\partial{\partial v^\alpha}
  \qquad
  \P^\nabla_{(x,v)}\Big(\frac\partial{\partial x^\mu}\Big)
  \;\;=\;\;
  \sum_{\alpha\,=\,1}^n\Gamma^\alpha_\mu(\,x,\,v\,)
  \,\frac\partial{\partial v^\alpha}
 \end{equation}
 with local coefficient functions $\Gamma^\alpha_\mu(\,x,\,v\,)$ generalizing
 the classical Christoffel symbols, because the coordinate vector fields
 $\frac\partial{\partial v^1},\,\ldots,\,\frac\partial{\partial v^n}$ span the
 vertical tangent bundle in every point of $\pi^{-1}(\,U\,)\,\subseteq\,VM$.\\

 With respect to the local trivialization of the vector bundle $VM$ provided
 by the local basis $\mathbf{v}_1,\,\ldots,\,\mathbf{v}_n$ every section
 $v\,\in\,\Gamma(\,M,\,VM\,)$ can be expanded as a sum $v\,=\,v^1\mathbf{v}_1
 \,+\,\ldots\,+\,v^n\mathbf{v}_n$ with local coefficient functions
 $v^1,\,\ldots,\,v^n$ in $C^\infty(U)$, which of course reappear in the
 local coordinate description of the section $v\,\in\,\Gamma(\,M,\,VM\,)$
 considered as a smooth map:
 $$
  (\,x^1,\,\ldots,\,x^m\,)\;\longmapsto\;(\,x^1,\,\ldots,\,x^m;
  \,v^1(x^1,\ldots,x^m),\,\ldots,\,v^n(x^1,\ldots,x^m)\,)\ .
 $$
 
 Calculating the differential of $v$ in these coordinates is straightforward:
 \begin{eqnarray*}
  \P^\nabla\left(v_*\,\Big(\,\frac\partial{\partial x^\mu}\,\Big)\right)
  &=&
  \P^\nabla\Big(\;\frac\partial{\partial x^\mu}\;+\;
  \sum_{\alpha\,=\,1}^n\frac{\partial v^\alpha}{\partial x^\mu}(\,x\,)
  \,\frac\partial{\partial v^\alpha}\;\Big)
  \\
  &=&
  \sum_{\alpha\,=\,1}^n
  \Big(\;\frac{\partial v^\alpha}{\partial x^\mu}(\,x\,)\;+\;
  \Gamma_\mu^\alpha(\,x,\,v(x)\,)\;\Big)\;\frac\partial{\partial v^\alpha}\ .
 \end{eqnarray*}
 
\noindent  Pushing this expression to a section of $VM$ via the vertical trivialization
 $\vtriv$ in diagram (\ref{linco}) converts the coordinate vector field
 $\frac\partial{\partial v^\alpha}$ back to the corresponding local section
 $\mathbf{v}^\alpha$, hence it reduces the covariant derivative to:
 \begin{equation}\label{cex}
  \nabla_{\frac\partial{\partial x^\mu}}v
  \;\;:=\;\;
  \sum_{\alpha\,=\,1}^n
  \Big(\;\frac{\partial v^\alpha}{\partial x^\mu}(\,x\,)\;+\;
  \Gamma_\mu^\alpha(\,x,\,v(x)\,)\;\Big)\;\mathbf{v}_\alpha\ .
 \end{equation}
 Directly from this expansion we conclude that a connection $\P^\nabla$ on
 $VM$ will not be a linear connection in the sense of Definition \ref{linx}
 unless its Christoffel symbols $\Gamma^\beta_\mu(\,x,\,v\,)$ are linear in
 the coefficient coordinates with an expansion
 \begin{equation}\label{lineq}
  \Gamma^\alpha_\mu(\,x^1,\,\ldots,\,x^m;\,v^1,\,\ldots,\,v^n\,)
  \;\;=\;\;
  \sum_{\omega\,=\,1}^n
  \Gamma^\alpha_{\mu\omega}(\,x^1,\,\ldots,\,x^m\,)\,v^\omega
 \end{equation}
 and suitable local coefficient functions $\Gamma^\alpha_{\mu\omega}(\,x\,)$.
 In the chosen linear coordinates the scalar multiplication $\Lambda_\lambda:\,
 VM\longrightarrow VM,\,v\longmapsto \lambda v,$ with $\lambda\,\in\,\R$ reads
 $(\,x^1,\ldots,x^m;\,v^1,\ldots,v^n\,)\longmapsto(\,x^1,\ldots,x^m;\,
 \lambda v^1,\ldots,\lambda v^n\,)$. It is thus a parallel endomorphism of
 $VM$ with respect to the connection $\P^\nabla$, if and only if we have
 for all $\mu$ equality between the two tangent vector expressions
 \begin{eqnarray*}
  \lefteqn{(\,\Lambda_\lambda\,)_{*,\,(x,v)}
  \Big(\,\frac\partial{\partial x^\mu}-\sum_{\alpha\,=\,1}^n
  \Gamma^\alpha_\mu(\,x,\,v\,)\frac\partial{\partial v^\alpha}\,\Big)}
  \quad
  &&
  \\
  &=&
  \frac\partial{\partial x^\mu}-\sum_{\alpha\,=\,1}^n
  \lambda\;\Gamma^\alpha_\mu(\,x,\,v\,)\frac\partial{\partial v^\alpha}
  \;\;\stackrel?=\;\;
  \frac\partial{\partial x^\mu}-\sum_{\alpha\,=\,1}^n
  \Gamma^\alpha_\mu(\,x,\,\lambda v\,)\frac\partial{\partial v^\alpha}
 \end{eqnarray*}
 due to $(\,\Lambda_\lambda\,)_*\frac\partial{\partial x^\mu}
 \,=\,\frac\partial{\partial x^\mu}$ and $(\,\Lambda_\lambda\,)_*
 \frac\partial{\partial v^\alpha}\,=\,\lambda\,\frac\partial{\partial
 v^\alpha}$. Note that the tangent vectors on the left and right hand
 side of this equation are horizontal in the points $(\,x,\,v\,)$ and
 $(\,x,\,\lambda v\,)$ of $VM$ respectively.\\

 In consequence the scalar multiplications $\Lambda_\lambda$ are parallel
 endomorphisms of $VM$ for all scalars $\lambda\,\in\,\R$, if and only if
 all the generalized Christoffel symbols $\Gamma^\alpha_\mu(x,v)$ are
 homogeneous functions of degree one in the coefficient coordinates
 $v^1,\,\ldots,\,v^n$. On the other hand every smooth homogeneous function
 of degree one is actually a linear function, hence we have proved:
 
 \begin{Lemma}[Characterization of Linear Connections]
 \hfill\label{clcx}\break
  A non--linear connection $\P^\nabla$ on a vector bundle gives rise
  to an $\R$--bilinear covariant derivative $\nabla:\,\Gamma(\,M,\,TM\,)
  \times\Gamma(\,M,\,VM\,)\longrightarrow\Gamma(\,M,\,VM\,)$, if and only
  if the multiplication by every $\lambda\,\in\,\R$ is a parallel endomorphism:
  $$
   \Lambda_\lambda:\;\;VM\;\longrightarrow\;VM,
   \qquad v\;\longmapsto\;\lambda\,v\ .
  $$
 \end{Lemma}
 
\begin{appendix}
\section{Explicit Calculations in Local Coordinates}
\label{excalc}
 In order to illustrate the rather abstract and functorial discussion of
 non--linear connections and their curvature in this article we will
 repeat parts of the arguments of Section $2$ in this appendix
 by explicit calculations in local coordinates. In particular we will
 use these calculations in local coordinates to prove the general curvature
 identity (\ref{cross}) relating the Nijenhuis tensor de\-finition (\ref{rcurv})
 of the curvature $R^\nabla$ of a non--linear connection on a fiber bundle
 $\F M$ to the commutator of iterated covariant derivatives of its sections.

 \pfill
 For actual calculations the most convenient definition of the tangent bundle
 of a smooth manifold $\F$ is as a set $T\F\,:=\,\mathrm{Jet}^1_0(\,\R,\,
 \F\,)$ of equivalence classes of curves $\gamma:\,\R\longrightarrow\F$
 under first order contact in $0\,\in\,\R$. The second iterated tangent
 bundle can be defined analogously as the set $T(\,T\F\,)$ of equivalence
 classes of maps $\gamma:\,\R^2\longrightarrow\F,\,(t,\varepsilon)\longmapsto
 \gamma(t,\varepsilon),$ under a modified equivalence relation of second order
 contact in $(0,0)$: Two maps $\gamma:\,\R^2\longrightarrow\F$ and
 $\hat\gamma:\,\R^2\longrightarrow\F$ determine the same equivalence class
 \begin{equation}\label{tc}
  \left.\frac d{dt}\right|_0\left.\frac d{d\varepsilon}\right|_0
  \gamma(\,t,\,\varepsilon\,)
  \;\;=\;\;
  \left.\frac d{dt}\right|_0\left.\frac d{d\varepsilon}\right|_0
  \hat\gamma(\,t,\,\varepsilon\,)
  \;\;\in\;\;
  T(\,T\F\,)\ ,
 \end{equation}
 if and only if for some and hence every local coordinate chart $(\,f^1,\,
 \ldots,\,f^n\,)$ of $\F$ defined in a neighborhood of $\gamma(\,0,0\,)
 \,=\,\hat\gamma(\,0,0\,)$ the $4n$ real numbers 
 \begin{equation}\label{2ct}
  \begin{array}{rcl}
   f^\alpha\Big(\;{\displaystyle\left.\frac d{dt}\right|_0\left.
   \frac d{d\varepsilon}\right|_0}\gamma(\,t,\,\varepsilon\,)\;\Big)
   &=&
   \quad(\,f^\alpha\circ\gamma\,)(0,0)
   \\[10pt]
   \dot f^\alpha\Big(\;{\displaystyle\left.\frac d{dt}\right|_0\left.
   \frac d{d\varepsilon}\right|_0}\gamma(\,t,\,\varepsilon\,)\;\Big)
   &=&
   \;\;{\displaystyle\left.\frac\partial{\partial t}\;\;\right|_{(0,0)}}
   (\,f^\alpha\circ\gamma\,)  
   \\[10pt]
   \vary f^\alpha\Big(\;{\displaystyle\left.\frac d{dt}\right|_0\left.
   \frac d{d\varepsilon}\right|_0}\gamma(\,t,\,\varepsilon\,)\;\Big)
   &=&
   \;\;{\displaystyle\left.\frac\partial{\partial\varepsilon}\;\;
   \right|_{(0,0)}}(\,f^\alpha\circ\gamma\,)  
   \\[10pt]
   \vary\dot f^\alpha\Big(\;{\displaystyle\left.\frac d{dt}\right|_0\left.
   \frac d{d\varepsilon}\right|_0}\gamma(\,t,\,\varepsilon\,)\;\Big)
   &=&
   {\displaystyle\left.\frac{\partial^2}{\partial t\,\partial\varepsilon}
   \right|_{(0,0)}}(\,f^\alpha\circ\gamma\,)
  \end{array}
 \end{equation}
 for $\gamma$ agree with the corresponding numbers for $\hat\gamma$. In
 this case these $4n$ numbers are of course the coordinates of the common
 equivalence class (\ref{tc}) in the induced coordinate chart $(\,f^\alpha;
 \,\dot f^\alpha;\,\vary f^\alpha;\,\vary\dot f^\alpha\,)$ of $T(\,T\F\,)$.
 In contrast to $T\F$ the manifold $T(\,T\F\,)$ is not a true jet space, the
 modification (\ref{2ct}) of the second order contact relation depends on
 the existence of tautological coordinates $(\,t,\varepsilon\,)$ on $\R^2$.
 From the definition (\ref{2ct}) of the equivalence relation we see directly
 that the involution $\Theta:\,T(T\F)\longrightarrow T(T\F)$ defined in
 equation (\ref{theta}) by precomposing $\gamma:\,\R^2\longrightarrow\F$
 with the swap $(t,\,\varepsilon)\longmapsto(\varepsilon,\,t)$ is
 well--defined on equivalence classes and takes the simple form
 \begin{equation}\label{thetaco}
  \Theta\,
  (\;f^\alpha;\,\dot f^\alpha;\,\vary f^\alpha;\,\vary{\dot f}^\alpha\;)
  \;\;=\;\;
  (\;f^\alpha;\,\vary f^\alpha;\,\dot f^\alpha;\,\vary{\dot f}^\alpha\;)
 \end{equation}
 in local coordinates, compare equation (\ref{picks}). Similarly the double
 projection $\Pi:\,T(\,T\F\,)\longrightarrow T\F\oplus T\F$ becomes in the
 chosen local coordinates
 $$
  \Pi(\;f^\alpha;\,\dot f^\alpha;\,\vary f^\alpha;\,\vary{\dot f}^\alpha\;)
  \;\;=\;\;
  (\;f^\alpha;\,\vary f^\alpha\;)\oplus(\;f^\alpha;\,\dot f^\alpha\;)
  \;\;=\;\;
  (\;f^\alpha;\,\vary f^\alpha;\,\dot f^\alpha\;)\ ,
 $$
 because $\oplus$ for vector bundles over $\F$ is essentially the same as
 $\times_\F$ and effectively removes the redundant copy of the coordinates
 $f^\alpha$. For an alternative coordinate chart $(\,h^1,\,\ldots,\,h^n\,)$
 of the manifold $\F$ the chain rule applied to the composition $h^\omega
 \circ\gamma\,=\,(h^\omega\circ f^{-1})\circ(f\circ\gamma)$ tells us
 $$
  \vary\dot h^\omega(\;f;\,\dot f;\,\vary f;\,\vary\dot f\;)
  \;\;=\;\;
  \sum_{\alpha,\,\beta\,=\,1}^n
  \frac{\partial^2h^\omega}{\partial f^\alpha\,\partial f^\beta}(\,f\,)
  \;\vary f^\alpha\;\dot f^\beta
  \;+\;
  \sum_{\alpha\,=\,1}^n\frac{\partial h^\omega}{\partial f^\alpha}(\,f\,)
  \;\vary\dot f^\alpha
 $$
 for the induced change of coordinates on $T(\,T\F\,)$ with similar
 formulas for $\dot h^\omega$ and $\vary h^\omega$. It is the symmetry
 of $\frac{\partial^2h^\omega}{\partial f^\alpha\,\partial f^\beta}(\,f\,)$
 under the interchange $\alpha\,\leftrightarrow\,\beta$, which allows
 $\Theta$ to read (\ref{thetaco}) in all local coordinates on $T(\,T\F\,)$
 induced by coordinates on $\F$! In taking differences the quadratic term
 drops out, hence the fibers of $\Pi$ are affine spaces modelled on $T\F$
 under the difference map
 \begin{equation}\label{diff}
  (\;f^\alpha;\,\dot f^\alpha;\,\vary f^\alpha;\;\mathbb{F}^\alpha\;)
  \;-\;
  (\;f^\alpha;\,\dot f^\alpha;\,\vary f^\alpha;\;\mathbb{\hat F}^\alpha\;)
  \;\;:=\;\;
  (\;f^\alpha;\,\mathbb{F}^\alpha\,-\,\mathbb{\hat F}^\alpha\;)
 \end{equation}
 in all local coordinates on $T(\,T\F\,)$ induced by local coordinates on $\F$.

 \pfill
 Turning to fiber bundles over a manifold $M$ we choose local coordinates
 on the total space of a fiber bundle by lifting local coordinates
 $(\,x^1,\,\ldots,\,x^m\,)$ of $M$ to functions on $\F M$ still denoted
 by $x^1,\,\ldots,\,x^m$. Complemented by suitable functions $f^1,\,\ldots,
 \,f^n\,\in\,C^\infty_\loc(\,\F M\,)$ we obtain a system of local coordinates
 $(\,x^1,\,\ldots,\,x^m;\,f^1,\,\ldots,\,f^n\,)$ with the property that the
 differentials $dx^1,\,\ldots,\,dx^m$ of the coordinates constant along the
 fibers of $\F M$ generate the horizontal forms so that the vertical tangent
 bundle is spanned by:
 $$
  \Vert_{(\,x,\,f\,)}\F M
  \;\;:=\;\;
  \mathrm{span}_\R\Big\{\;
  \frac\partial{\partial f^1},\,\ldots,\,\frac\partial{\partial f^n}\;\Big\}
  \;\;\subseteq\;\;
  T_{(\,x,\,f\,)}\F M
 $$
 All the constructions discussed in the first part of this appendix for the
 second iterated tangent bundle $T(\,T\F\,)$ extend without further ado to
 the second iterated vertical tangent bundle $\Vert\,\Vert\,\F M$ by throwing
 in the additional coordinates $(\,f^\alpha;\,\dot f^\alpha;\,\vary f^\alpha;\,
 \vary\dot f^\alpha\,)\,\rightsquigarrow\,(\,x^\mu;\,f^\alpha;\,\dot f^\alpha;
 \,\vary f^\alpha;\,\vary\dot f^\alpha\,)$ of the base manifold $M$ as passive
 coordinates. In this specific local coordinate setup for a fiber bundle an
 arbitrary non--linear connection $\P^\nabla$ reads
 $$
  \P^\nabla_{(\,x,f\,)}\Big(\frac\partial{\partial f^\alpha}\Big)
  \;\;=\;\;
  \frac\partial{\partial f^\alpha}
  \qquad\quad
  \P^\nabla_{(\,x,f\,)}\Big(\frac\partial{\partial x^\mu}\Big)
  \;\;=\;\;
  \sum_{\alpha\,=\,1}^n\,\Gamma_\mu^\alpha
  (\,x,\,f\,)\;\frac\partial{\partial f^\alpha}
 $$
 with suitable local coefficient functions $\Gamma^\alpha_\mu(\,x,\,f\,)$
 generalizing the standard Christoffel symbols. Needless to say the curvature
 $R^\nabla$ of the non--linear connection $\P^\nabla$ is determined by these
 generalized Christoffel symbols, due to horizontality, it suffices to calculate
 expression (\ref{rcurv}) for the vector fields:
 \begin{eqnarray*}
  \lefteqn{R^\nabla\Big(\,\frac\partial{\partial x^\mu},\,
  \frac\partial{\partial x^\nu}\,\Big)}
  &&
  \\
  &=&
  -\,\P^\nabla\left[\;
  \frac\partial{\partial x^\mu}\;-\;\sum_{\alpha\,=\,1}^n
  \Gamma^\alpha_\mu(\,x,f\,)\,\frac\partial{\partial f^\alpha},\;
  \frac\partial{\partial x^\nu}\;-\;\sum_{\beta\,=\,1}^n
  \Gamma^\beta_\nu(\,x,f\,)\,\frac\partial{\partial f^\beta}\;\right]
  \\
  &=&
  +\;\sum_{\beta\,=\,1}^n
  \frac{\partial\Gamma^\beta_\nu}{\partial x^\mu}(\,x,f\,)
  \,\frac\partial{\partial f^\beta}
  \;-\;
  \sum_{\alpha,\,\beta\,=\,1}^n\Gamma^\alpha_\mu(\,x,f\,)\,
  \frac{\partial\Gamma^\beta_\nu}{\partial f^\alpha}(\,x,f\,)\,
  \frac\partial{\partial f^\beta}
  \\
  &&
  -\;\sum_{\alpha\,=\,1}^n
  \frac{\partial\Gamma^\alpha_\mu}{\partial x^\nu}(\,x,f\,)
  \,\frac\partial{\partial f^\alpha}\;+\;
  \sum_{\beta,\,\alpha\,=\,1}^n\Gamma^\beta_\nu(\,x,f\,)\,
  \frac{\partial\Gamma^\alpha_\mu}{\partial f^\beta}(\,x,f\,)\,
  \frac\partial{\partial f^\alpha}\ .
 \end{eqnarray*}
 Interchanging the dummy indices $\alpha\leftrightarrow\beta$ if necessary,
 this results in a sum of the form $R^\nabla(\,\frac\partial{\partial x^\mu},
 \,\frac\partial{\partial x^\nu}\,)\,=\,\sum R^\alpha_{\mu\nu}(\,x,\,f\,)\,
 \frac\partial{\partial f^\alpha}$ with coefficients given by:
 \begin{equation}\label{cloc}
  R^\alpha_{\mu\nu}
  \;\;:=\;\;
  \frac{\partial\Gamma^\alpha_\nu}{\partial x^\mu}
  \;-\;\frac{\partial\Gamma^\alpha_\mu}{\partial x^\nu}
  \;+\;\sum_{\beta\,=\,1}^n\Big(\;
  \Gamma^\beta_\nu\,\frac{\partial\Gamma^\alpha_\mu}{\partial f^\beta}
  \;-\;
  \Gamma^\beta_\mu\,\frac{\partial\Gamma^\alpha_\nu}{\partial f^\beta}
  \;\Big)
 \end{equation}
 In case the fiber bundle is actually a vector bundle we have a distinguished
 class of local coordinates $(\,x^1,\ldots,x^m;\,v^1,\ldots,v^n\,)$ on the
 total space $VM$ associated to a basis $\mathbf{v}_1,\,\ldots,\,\mathbf{v}_n
 \,\in\,\Gamma_\loc(\,M,VM\,)$ of local sections:
 $$
  (\,x^1,\,\ldots\,x^m;\,v^1,\,\ldots,\,v^n\,)
  \;\;\longleftrightarrow\;\;
  v^1\,\mathbf{v}_1\;+\;\ldots\;+\;v^n\,\mathbf{v}_n
  \;\;\in\;\;
  V_{(x^1,\ldots,x^m)}M\ .
 $$
 In these linear local coordinates on $VM$ the generalized Christoffel
 symbols $\Gamma^\alpha_\mu(\,x,\,v\,)$ of a linear connection are
 linear functions with an expansion
 $$
  \Gamma^\alpha_\mu(\,x^1,\,\ldots,\,x^m;\,v^1,\,\ldots,\,v^n\,)
  \;\;=\;\;
  \sum_{\omega\,=\,1}^n
  \Gamma^\alpha_{\mu\omega}(\,x^1,\,\ldots,\,x^m\,)\,v^\omega
 $$
 with local coefficient functions $\Gamma^\alpha_{\mu\omega}(\,x\,)$. In
 particular the expansion (\ref{cloc}) of the curvature of a linear
 connection simplifies in linear coordinates on $VM$
 $$
  R^\nabla\Big(\;\frac\partial{\partial x^\mu},
  \;\frac\partial{\partial x^\nu}\;\Big)
  \;\;=\;\;
  \sum_{\alpha,\,\omega\,=\,1}^nR_{\mu\nu;\,\omega}^\alpha
  (\,x^1,\,\ldots,\,x^m\,)\;v^\omega\;\frac\partial{\partial v^\alpha}
 $$
 with coefficients $R_{\mu\nu;\,\omega}^\alpha(\,x\,)$ given by the
 classical formula:
 \begin{equation}\label{ccurv}
  R_{\mu\nu;\,\omega}^\alpha
  \;\;=\;\;
  \frac{\partial\Gamma_{\nu\omega}^\alpha}{\partial x^\mu}
  \;-\;\frac{\partial\Gamma_{\mu\omega}^\alpha}{\partial x^\nu}
  \;+\;\sum_{\beta\,=\,1}^n
  \Big(\;\Gamma_{\mu\beta}^\alpha\,\Gamma_{\nu\omega}^\beta
  \;-\;\Gamma_{\nu\beta}^\alpha\,\Gamma_{\mu\omega}^\beta\;\Big)\ .
 \end{equation}
 Coming back to fiber bundles we want to discuss the covariant
 derivative of a local section $f\,\in\,\Gamma_\loc(\,M,\F M\,)$ of
 $\F M$, which in the chosen local coordinates on $M$ and $\F M$ becomes
 an $n$--tuple of smooth functions:
 $$
  f(\,x^1,\ldots,x^m\,)
  \;\;=\;\;
  (\,x^1,\ldots,x^m;\,f^1(x^1,\ldots,x^m),\ldots,f^n(x^1,\ldots,x^m)\,)\ .
 $$
 In turn its covariant derivative with respect to the connection $\P^\nabla$
 reads:
 \begin{eqnarray*}
  \Big(\,D^\nabla_{\frac\partial{\partial x^\mu}}f\,\Big)(\,x\,)
  &=&
  \P^\nabla_{(\,x,f(x)\,)}\Big(\;\frac\partial{\partial x^\mu}
  \;+\;\sum_{\alpha\,=\,1}^n\frac{\partial f^\alpha}{\partial x^\mu}(\,x\,)
  \frac\partial{\partial f^\alpha}\;\Big)
  \\
  &=&
  \sum_{\alpha\,=\,1}^n\Big(\;\frac{\partial f^\alpha}{\partial x^\mu}(\,x\,)\;
  +\;\Gamma^\alpha_\mu(\,x,f(x)\,)\;\Big)\,\frac\partial{\partial f^\alpha}\ ,
 \end{eqnarray*}
 In this formula we have considered the result as a vertical vector field,
 alternatively we may write the result as a local section of $\Vert\,\F M$
 \begin{equation}\label{gcx}
  \Big(\,D^\nabla_{\frac\partial{\partial x^\mu}}f\,\Big)(\,x\,)
  \;\;=\;\;
  \Big(\;x;\,f^\alpha(x);\,\frac{\partial f^\alpha}{\partial x^\mu}(x)
  \,+\,\Gamma^\alpha_\mu(\,x,\,f(x)\,)\;\Big)\ .
 \end{equation}
 For a smooth family $(\,f_\varepsilon\,)$ of local sections of $\F M$
 we obtain the variation
 \begin{eqnarray*}
  \Big({\textstyle\left.\frac d{d\varepsilon}\right|_0}
  D^\nabla_{\frac\partial{\partial x^\mu}}f_\varepsilon\Big)(\,x\,)
  &=&
  \Big(\;x;\,f^\alpha(x);\,\frac{\partial f^\alpha}{\partial x^\mu}(x)
  \,+\,\Gamma^\alpha_\mu(\,x,\,f(x)\,);
  \\
  &&
  \vary f^\alpha(x);\,\frac{\partial\,\vary f^\alpha}{\partial x^\mu}(x)
  +\sum_{\beta\,=\,1}^n\frac{\partial\Gamma^\alpha_\mu}{\partial f^\beta}
  (\,x,f(x)\,)\,\vary f^\beta(x)\;\Big)
 \end{eqnarray*}
 with $f^\alpha\,:=\,f^\alpha_0$ and $\vary f^\alpha\,:=\,\left.
 \frac d{d\varepsilon}\right|_0f^\alpha_\varepsilon$. Evidently we
 need to apply the fiberwise involution $\Theta$ to interpret the
 result as the covariant derivative of the variation $\left.\frac d
 {d\varepsilon}\right|_0f_\varepsilon\,\in\,\Gamma_\loc(\,M,\Vert\,\F M\,)$
 with respect to the non--linear connection $\P^{\nabla^\Vert}$ induced by
 equation (\ref{vconex}) on the vertical tangent bundle:
 \begin{eqnarray*}
  \Big(\;D^{\nabla^\Vert}_{\frac\partial{\partial x^\mu}}
   {\textstyle\left.\frac d{d\varepsilon}\right|_0f_\varepsilon}\;\Big)(\,x\,)
  &=&
  \Big(\;x;\,f^\alpha(x);\,\vary f^\alpha(x);\,\frac{\partial f^\alpha}
  {\partial x^\mu}(x)\,+\,\Gamma^\alpha_\mu(\,x,\,f(x)\,);
  \\
  &&
  \qquad\frac{\partial\,\vary f^\alpha}{\partial x^\mu}(x)
  \,+\,\sum_{\beta\,=\,1}^n\frac{\partial\Gamma^\alpha_\mu}{\partial f^\beta}
  (\,x,f(x)\,)\,\vary f^\beta(x)\;\Big).
 \end{eqnarray*}
 Comparing the resulting expression with the formula (\ref{gcx}) for a general
 non--linear connection we conclude that the connection $\P^{\nabla^\Vert}$ is
 given by
 $$
  \P^{\nabla^\Vert}_{(\,x,f,\vary f\,)}
  \Big(\,\frac\partial{\partial x^\mu}\,\Big)
  \;\;=\;\;
  \sum_{\alpha\,=\,1}^n\,\Gamma_\mu^\alpha
  (\,x,f\,)\;\frac\partial{\partial f^\alpha}
  \;+\;
  \sum_{\alpha,\,\beta\,=\,1}^n
  \,\frac{\partial\Gamma_\mu^\alpha}{\partial f^\beta}
  (\,x,f\,)\,\vary f^\beta\,\frac\partial{\partial\,\vary f^\alpha}\ ,
 $$
 note that the vertical tangent vectors $\frac\partial{\partial f^\alpha}$ and
 $\frac\partial{\partial\,\vary f^\alpha}$ are fixed by every non--linear
 connection on $\Vert\,\F M$. Substituting the variation $\vary f^\alpha(x)$
 by the corresponding expression $\frac{\partial f^\alpha}{\partial x^\nu}(x)
 \,+\,\Gamma^\alpha_\nu(x,f(x))$ in the covariant derivative
 $D^\nabla_{\frac\partial{\partial x^\nu}}f$ of the local section
 $f\,\in\,\Gamma_\loc(\,M,\F M\,)$ we obtain directly:
 \begin{eqnarray*}
  \lefteqn{\Big(\,D^{\nabla^\Vert}_{\frac\partial{\partial x^\mu}}
  D^{\nabla^{\hphantom{V}}}_{\frac\partial{\partial x^\nu}}f\,\Big)(\,x\,)}
  \quad
  &&
  \\
  &=&
  \Big(\;x;\,f^\alpha(x);\,\frac{\partial f^\alpha}{\partial x^\nu}(x)
  \,+\,\Gamma^\alpha_\nu(\,x,f(x)\,);\,\frac{\partial f^\alpha}
  {\partial x^\mu}(x)\,+\,\Gamma^\alpha_\mu(\,x,f(x)\,);
  \\[2pt]
  &&
  \;\;+\,\frac{\partial^2f^\alpha}{\partial x^\mu\partial x^\nu}(x)
  \,+\,\frac{\partial\Gamma^\alpha_\nu}{\partial x^\mu}(x)\,+\,
  \sum_{\beta\,=\,1}^n\frac{\partial\Gamma^\alpha_\mu}{\partial f^\beta}
  (\,x,f(x)\,)\,\Gamma^\beta_\nu(\,x,f(x)\,)
  \\
  &&
  \;\;+\,\sum_{\beta\,=\,1}^n\frac{\partial\Gamma^\alpha_\nu}
  {\partial f^\beta}(\,x,f(x)\,)\,\frac{\partial f^\beta}{\partial x^\mu}(x)
  \,+\,\sum_{\beta\,=\,1}^n\frac{\partial\Gamma^\alpha_\mu}{\partial f^\beta}
  (\,x,f(x)\,)\,\frac{\partial f^\beta}{\partial x^\nu}(x)\;\Big)\ .
 \end{eqnarray*}
 Evidently it is necessary to apply the involution $\Theta$ first to form the
 difference
 \begin{eqnarray*}
  \lefteqn{\left(\;D^{\nabla^\Vert}_{\frac\partial{\partial x^\mu}}
  D^{\nabla^{\hphantom{V}}}_{\frac\partial{\partial x^\nu}}f
  \;-\;\Theta\Big(\,D^{\nabla^\Vert}_{\frac\partial{\partial x^\nu}}
  D^{\nabla^{\hphantom{V}}}_{\frac\partial{\partial x^\mu}}f\,\Big)
  \;\right)(\,x\,)}
  \qquad
  \\[2pt]
  &=&
  \Big(\;x;\,f^\alpha(x);\,
  \frac{\partial\Gamma^\alpha_\nu}{\partial x^\mu}
  \,-\,\frac{\partial\Gamma^\alpha_\mu}{\partial x^\nu}
  \,+\,\sum_{\beta\,=\,1}^n\Big(\,
  \Gamma^\beta_\nu\,\frac{\partial\Gamma^\alpha_\mu}{\partial f^\beta}\,-\,
  \Gamma^\beta_\mu\,\frac{\partial\Gamma^\alpha_\nu}{\partial f^\beta}\,\Big)
  \;\Big)
  \\
  &=&
  \Big(\;x;\,f^\alpha(x);\,R^\alpha_{\mu\nu}(\,x,\,f(x)\,)\;\Big)
  \;\;=\;\;
  \Big(\;R^\nabla_{\frac\partial{\partial x^\mu},
   \,\frac\partial{\partial x^\nu}}f\;\Big)(\,x\,)\ ,
 \end{eqnarray*}
 where the coefficients $R^\alpha_{\mu\nu}(x,f)$ of the curvature of the
 connection $\P^\nabla$ have been calculated in equation (\ref{cloc}). Due
 to $[\,\frac\partial{\partial x^\mu},\,\frac\partial{\partial x^\nu}\,]\,=\,0$
 the latter identity is exactly the local coordinates formulation of the
 curvature identity (\ref{cross}).

\end{appendix}

\end{document}